\documentclass[preprint,12pt]{amsart}
\usepackage{geometry}                
\usepackage{graphicx}
\usepackage{amssymb}
\usepackage{epstopdf}
\usepackage{amsthm}
\usepackage{amsmath}

\newcommand{\ind}{1\!\!1}

\newcommand{\K}{\mathcal{K}}
\newcommand{\onto}[2]{\genfrac{}{}{0pt}{}{#1}{#2}}
\newtheorem*{hyp}{Assumption}
\newtheorem{lem}{Lemma}

\newtheorem{thm}{Theorem}

\title{A Kriging procedure for processes indexed by graphs}
\author[]{T. Espinasse \& J-M. Loubes \\ $ $ \\ Institut de Math\'ematiques de Toulouse, France}

\begin{document}
\maketitle
\begin{abstract}
We provide a new  kriging procedure of processes on graphs. Based on the construction of Gaussian random processes indexed by graphs, we extend to this framework the usual linear prediction method for spatial random fields, known as kriging. We provide the expression of the estimator of such a random field at unobserved locations as well as a control for the prediction error.
\end{abstract}


\section*{Introduction}

Data presenting an inner geometry, among them data indexed on graphs are growing tremendously in size and prevalence these days. High dimensional  data often present  structural links that can be modeled  through a graph representation, for instance the World Wide Web graph or the  social networks well studied in history or geography \cite{moyenage}, or molecular graphs in biology or medicine for instance. The definition and the analysis of  processes indexed by such graphs is of growing interest in the statistical community. In particular, the definition of graphical models by J.N. Darroch, S.L. Lauritzen and T.P. Speed in $1980$ \cite{graph_mod} fostered new interest in Markov fields, and many tools have been developed in this direction (see, for instance~\cite{verzelen1} and~\cite{verzelen2}). When confronted to missing data or to forecast the values of the process at unobserved sites, many prediction methods have been proposed in the statistics community over the past few years. Among them,  we focus in 
this paper on a new method extending the kriging procedure to the case of a Gaussian process indexed by a graph.\vskip .1in 

Kriging is named for the mining engineer Krige, whose paper \cite{Krige51} introduced the method. For background on kriging see \cite{Ste} or \cite{MR1239641}. Gaussian process models Ðalso called Kriging modelsÐ are often used as mathematical approximations of expensive experiments. Originally presented in spatial statistics  as an optimal linear unbiased predictor of random processes, its interpretation is usually restricted to the convenient framework of Gaussian Processes (GP). It is based on the computation of the conditional expectancy, which requires a proper definition of the covariance structure of the process. For this, we will use the spectral definition of a covariance operator based on the adjacency operator. Inspired by \cite{2011arXiv1104.3664E}, we define Gaussian fields on graphs using their  spectral representation  (see for instance in \cite{AzDa}). Then, we extend to this case the usual method of prediction using the Kriging method. For this we will use blind prediction technics 
generalizing to graphs process the procedure for time series presented in~\cite{MMStibo}.  \vskip .1in

The paper falls into the following parts. Section~\ref{s:notations} is devoted to the definitions of a graph and the presentation of the problem. The forecast problem is tackled in Section~\ref{s:not blind}. The proofs are postponed to the Appendix.

\section{Preliminary notations and definitions} \label{s:notations}

\subsection{Random processes indexed by graphs}
In this section, we introduce both the context and the objects considered in the whole paper.

First, assume that $\mathbf{G}$ is a graph, that is a set of vertices $G$ and a set of edges $E \subset G \times G$. In this work, $G$ is assumed to be infinite (but countable).

Two vertices $i,j \in G$ are neighbors if $(i,j) \in E$. The degree $d(i)$ of a vertex $i \in G$ is the number of neighbors and the degree of the graph $\mathbf{G}$ is defined as the maximum degree of the vertices of the graph $\mathbf{G}$ :
$$\operatorname{deg}(\mathbf{G}) := \max_{i \in G} \operatorname{deg}(i).$$ 
From now on, we assume that the degree of the graph $\mathbf{G}$ is bounded, that is 
$$\exists d_{\text{max}} >0, \forall i \in G, d(i) \leq d_{\text{max}}. $$

Furthermore, the graph $\mathbf{G}$ is endowed with the natural distance $d_G$, that is the length of the shortest path between two vertices.\vskip .1in
In the following, we will consider  the renormalized adjacency operator $A$ of $\mathbf{G}$. So its entries belong to $[-\frac{1}{\operatorname{deg}(\mathbf{G})},\frac{1}{\operatorname{deg}(\mathbf{G})}]$. It is defined as
$$A_{ij} = \frac{1}{d_{\text{max}}} \ind_{i=j}. $$
We denote by $B_G$ the set of all bounded Hilbertian
 operators on $l^2(G)$ (the set of square sommable real sequences indexed by $G$). 

To introduce the spectral decomposition, consider the action of the adjacency operator on $l^2(G)$ as 
$$\forall u \in l^2(G), (Au)_i := \sum_{j \in G} A_{ij} u_j, (i\in G).$$

The operator space $B_G$ will be endowed with the classical operator norm 
$$\forall T \in B_G, \left\| T \right\|_{2,op}: = \sup_{u \in l^2(G), \left\| u\right\|_2 \leq 1} \left\| Tu \right\|_2 , $$
where $\left\| . \right\|_2$ stands for the usual norm on $l^2(G)$.

Notice that, as the degree of $\mathbf{G}$ and the entries of $A$ are both bounded, $A$ lies in $B_{G}$, and we have $$\left\| A \right\|_{2,op} \leq 1 .$$  
Finally we get that $A$ is  a symmetric bounded normal Hilbertian operator. Recall that for any bounded Hilbertian operator $A \in B_G$, the spectrum $\operatorname{Sp}(A)$ is defined as the set of all complex numbers $\lambda$ such that 
$\lambda \operatorname{Id}- A$ is not invertible (here $\operatorname{Id}$ stands for the identity on $l^2(G)$).  So the spectrum of the normalized adjacency operator $A$ is a non-empty compact subset of $\mathbb{R}$.


Using the spectral representation of the graph, we proved in \cite{2011arXiv1104.3664E} that we can define Gaussian processes indexed on graphs whose covariance structure relies only on the geometry of the graph via its spectrum. For this, for any bounded positive function $f$, analytic on the convex hull of $\operatorname{Sp}(A)$, note first that $f(A)$ defines a bounded positive definite symmetric operator on $l^2(G)$. 
Then we can define a Gaussian process $(X_i)_{i \in G}$ indexed by the vertices $G$ of the graph $\mathbf{G}$ with covariance operator $\Gamma$ defined as
$$\Gamma =\int_{\operatorname{Sp}(A)} f(\lambda) \mathrm{d} E (\lambda). $$ 
Such graph analytical processes extend the notion of time series to a graph indexed process. So using this terminology, we will say that $X$ is
\begin{itemize}
 \item $MA_q$ if $f$ is a polynomial of degree $q$.
\item $AR_p$ if $\frac{1}{f}$ is a polynomial of degree $p$ which has no root in the convex hull of $\operatorname{Sp}(A)$.
\item $ARMA_{p,q}$ if $f = \frac{P}{Q}$ with $P$ a polynomial of degree $p$ and $Q$ a polynomial of degree $q$ with no roots in the convex hull of $\operatorname{Sp}(A)$.
\end{itemize}
Otherwise, we will talk about the $MA_\infty$ representation of the process $\mathbf{X}$.
We call $f$ the {\it spectral density} of the process $\mathbf{X}$, and denote its corresponding covariance operator by 
$$\Gamma = \mathcal{K}(f)=f(A) $$

So the spectral analysis f the graph enables to  define a class of admissible covariances for stationary Gaussian processes with associated spectral density $f$. Hereafter we tackle the issue of prediction of such processes.

\subsection{Blind prediction problem}

We can now introduce our prediction problem on this framework.

The problem comes from practical issues. In real life problems, it is usual to own a single sample, which has to be used for both estimation and prediction.

Let  $f$ be a bounded positive function, analytic over $\operatorname{Sp}(A)$, and $(X_i)_{i\in G}$ be a Gaussian zero-mean process indexed by $G$ of covariance operator $\K(f)$.

We will observe the process on a growing sequence of subgraphs of $G$, but with missing values we aim at predicting. Let $(O,B)$ be a partition of $G$. The set $O$ will denote the set of indexes for all the possibly observed values while $B$ denotes the "blind" missing values index set. In all the following, the set $B$ where the observations should be forecast, is assumed to be finite.  

Let $(G_N)_{N \in \mathbb{N}}$ be a growing sequence of induced subgraphs of $G$. This means that we have at hand a growing sequence of vertices and consider as the observed graph all the existing edges between the vertices that are observed. From now on, we assume that $N$ is large enough to ensure $B \subset G_N$. The observation index set will be denoted $O_N := O \cap G_N$ . \\
Hence, we consider the restriction $X_{O_N} := (X_i)_{i \in O_N}$, which stands for the data we have at hand at step $N$.  We consider the asymptotic framework where the observations $O_N$  fill the space between the blind part and the graph, in the sense that the distance between the blind locations $B$ and the non observed graph $G \setminus O_N$ increases, i.e  $$ m_N := \frac{1}{d_G(B,G\setminus O_N)} \longrightarrow 0,\quad {\rm when} \: N \rightarrow + \infty.$$
This corresponds to the natural case where the blind part of the graph becomes more and more surrounded by the observations without any gap. So the non observed locations of the process tend to be closer to observations of the process, which implies that forecasting the values at the blind locations become possible and relies on the rate at which such gap is filled, namely $m_N$.
\vskip .1in

 In the following, we will make the assumption that we can  dispose of a consistent estimation procedure $\hat{f}_N$ for $f$, such that there exists a decreasing sequence $r_N \rightarrow 0$, with is  a rate of convergence of the spectral density when $N$ goes to infinity, providing the controls on the following estimation errors 

\begin{hyp} {Preliminar consistent estimate for the spectral density}

\begin{itemize}
\item $  \mathbb{E}\left[\left\| \hat{f}_{N}-f \right\|_\infty^2 \right]^{\frac{1}{2}} \leq r_N.$
\item $  \mathbb{E}\left[ \left\| \hat{f}_{N}-f \right\|_\infty^4 \right]^{\frac{1}{4}} \leq r_N.$
\end{itemize}
\end{hyp}


Our aim is, observing only one realization $X_{O_N}$, to perform both estimation and prediction of any variable $Z_B$ 
 defined as a linear combination of the process taken at unobserved locations, i.e of the following form 
$$Z_B=a_B^TX_B, $$ 
with $a_B \in \mathbb{R}^B, \left\| a_B \right\|_2 = 1$.

Hereafter, we will write, for sake of simplicity, extracted operators like block matrices even if they are of infinite size
\begin{align*}
\forall U,V \subset G, \K_{UV}(f) & = \left( \K(f)_{ij}\right)_{i \in U, j \in V}
\\ \forall V \subset G, \K_{V}(f) & = \left( \K(f)_{ij}\right)_{i , j \in V}
\end{align*}

Recall (see for instance in \cite{Ste}) that the best linear predictor of $Z_B$ (which is also the best predictor in the Gaussian case) can be written as
$$\bar{Z}_B = P_{[X_{O_N}]}(f) Z_B := a_B^T \K_{B O_N}(f)\left(\K_{O_N}(f)\right)^{-1}X_{O_N}.$$
$\K_{B O_N}(f)$ is  the covariance between the observed process and the blind part, while   $\K_{O_N}(f)$ corresponds to the covariance of the process restricted to the observed data points. 
Note that this projection term is well defined. Indeed since $f$ is positive, $\K(f)$ is invertible, and therefore, $\K_{O_N}(f)$ is also invertible, as a principle minor. However, since $f$ is unknown we can not use this as a direct estimator. 

Then, remark that, asymptotically, in the sense $N\rightarrow + \infty$, we  may observe the process at all locations, $X_O$. We thus can  introduce the best linear prediction of $Z_B$ knowing  all the possible observations $X_O$ as
$$\tilde{Z}_B :=P_{[X_{O}]}(f) Z_B := a_B^T \K_{B O}(f)\left(\K_{O}(f)\right)^{-1}X_{O}. $$

Finally, the blind forecasting problem can be formulated as a two step procedure mixing the estimation of the projector operator and the prediction using the estimated projector. It can be thus decomposed as follows
\begin{itemize}
\item Estimation step: estimate $ P_{[X_{O_N}]}(f)$ by $ \hat{P}_{[X_{O_N}]}(f):= P_{[X_{O_N}]}(\hat{f})$.
\item Prediction step: build $\hat{Z}_B := P_{[X_{O_N}]}(\hat{f})Z_B$.
\end{itemize}

Therefore, it seems natural, in order to analyze  the forecast procedure, to consider  an upper bound on the risk defined by 

$$\mathbf{R}_{N} =  \sup_{ \onto{ Z_B   = a_B^TX_B}{ \left\|a_B\right\|_2 = 1}}  \mathbb{E} \left[      \left( Z_B - \hat{Z}_B \right)^2               \right]  ^{\frac{1}{2}}.$$

But we can see that this risk admits the following decomposition
$$\mathbf{R}_{N}  = \sup_{ \onto{ Z_B   = a_B^TX_B}{ \left\|a_B\right\|_2 = 1}}  \left( \mathbb{E} \left[      \left( Z_B - \tilde{Z}_B \right)^2   \right]^{\frac{1}{2}}  + \mathbb{E} \left[      \left( \tilde{Z}_B - \hat{Z}_B \right)^2               \right]^{\frac{1}{2}} \right). $$

Then, notice that the first term of this sum does not decrease to $0$ when $N$ goes to infinity. Actually it is an innovation type term. Therefore, we  consider the upper bound

\begin{align*}
\mathbf{R}_{N} &\leq  \sup_{ \onto{ Z_B   = a_B^TX_B}{ \left\|a_B\right\|_2 = 1}} \mathbb{E} \left[      \left( Z_B - \tilde{Z}_B \right)^2               \right] ^{\frac{1}{2}}  +\sup_{ \onto{ Z_B   = a_B^TX_B}{ \left\|a_B\right\|_2 = 1}} \mathbb{E} \left[      \left( \tilde{Z}_B - \bar{Z}_B \right)^2               \right]^{\frac{1}{2}} 
\\  &   \hspace{1cm} + \sup_{ \onto{ Z_B   = a_B^TX_B}{ \left\|a_B\right\|_2 = 1}}  \mathbb{E} \left[      \left( \bar{Z}_B - \hat{Z}_B \right)^2               \right] ^{\frac{1}{2}}  \\ & \leq  \sup_{ \onto{ Z_B   = a_B^TX_B}{ \left\|a_B\right\|_2 = 1}} \mathbb{E} \left[      \left( Z_B - \tilde{Z}_B \right)^2               \right] ^{\frac{1}{2}}   + \mathcal{R}_{N},
\end{align*}
where we have set
$$\mathcal{R}_{N}  :=  \sup_{ \onto{ Z_B   = a_B^TX_B}{ \left\|a_B\right\|_2 = 1}} \mathbb{E} \left[      \left( \tilde{Z}_B - \bar{Z}_B \right)^2               \right]^{\frac{1}{2}}   +\sup_{ \onto{ Z_B   = a_B^TX_B}{ \left\|a_B\right\|_2 = 1}} \mathbb{E} \left[      \left( \bar{Z}_B - \hat{Z}_B \right)^2               \right] ^{\frac{1}{2}} .$$
The first term does not depend on the estimation procedure, hence our main issue is to compute the rate of convergence towards $0$ of the two last terms of the previous sum. That is the reason why, $\mathcal{R}_{N}$ plays the role of the estimation risk, that will be controlled  in the whole paper.

\section{Prediction of a graph process with  independent observations} \label{s:not blind}

In this section, we assume that another sample $Y$ independent of $X$ and drawn with the same distribution, is available. We will use $Y$ to perform the estimation, and plug this estimation in order to predict $Z_{B}$.

More precisely, assume that for $p\geq 1$, $\hat{f}_N$ is build using the sample $Y$ observed on $O_N$ and that $X$ is observed on another subsample of the graph $O_p$. These nested collection of subgraphs fulfills the condition that $\frac{1}{m_p} = d_G\left(B, (G\setminus O_p)\right)$ goes to infinity. In this first part, we thus have two independent asymptotics. The first one with respect to $N$ controls how close the estimated spectral density will be from the true one while the second with respect to $p$ deals with the accuracy of the linear projector onto $O_p$.

We want to control the prediction of $Z_B$ using the estimation $\hat{f}_N$.
%
%
%
%
%
%
%
%
%

To give an upper bound, we need to assume regularity for the spectral density $f$. First, we assume that the process has short range memory through the following assumption:

\begin{hyp}[A1]
There exists $M >0$ such that
$$\forall t \in \operatorname{Sp}(A), m \leq f(t) \leq M.$$
\end{hyp}

Actually, we will assume that all the estimators should verify the same inequality, that is, for all $N \in \mathbb{N}$,
\begin{hyp}[A2]
$$\forall t \in \operatorname{Sp}(A), m \leq \hat{f}_N(t) \leq M.$$
\end{hyp}

We will need another regularity assumption on $f$ to control high-frequency behavior.

\begin{hyp}[A3]
The function $f(x) = \sum f_k x^k.$  is analytic on a the compact disk $\bar{D}(0,1)$ and verifies
$$\sum k |f_k| \leq M. $$
\end{hyp}

Note that this assumption ensures that $\frac{1}{f}$ is absolutely convergent in $1$ and that, if $\frac{1}{f}(x) = \sum (\frac{1}{f})_k x^k$, 
$$\sum k |(\frac{1}{f})_k| \leq \frac{M}{m^2}. $$

The following theorem provides an upper bound for the risk $\mathcal{R}_{N,p}$

\begin{thm} \label{th:indep}
Assume that two independent samples are available. Under Assumptions (A1), (A2) and (A3), the risk admits the following upper bound
$$\mathcal{R}_{N,p} \leq   \frac{\sqrt{M}(m+M)}{m^2}  r_N +  \frac{M}{m^4}  \left( \frac{M^{\frac{5}{2}}}{m} + M^{\frac{3}{2}} \right)   m_p .$$
\end{thm}

The prediction risk $\mathcal{R}_{N,p}$ is made of two separate terms which go to zero independently  when $p$ and $N$ increase. The first term is entirely governed by the accuracy of the estimation of the spectral density. The second term depends on the decay of  $\sum_{k \geq d_G\left(B, (G\setminus G_p)\right)} \left| (\frac{1}{f})_k \right| $. This sum is the remaining term of a convergent sum and thus vanishes when $p$ growths large.
Note that using Assumption (A3), we obtain a bound in $m_p$. Adding more regularity on the function $f$ of the type  $$\sum k^{2s} |(\frac{1}{f})_k| \leq \frac{M}{m^2}$$ for a given $s>0$   would improve this rate by replacing $m_p$ by $m_p^s$. \vskip .1in 

\begin{proof}
The proof consists in bounding the two terms in the decomposition of the risk $\mathcal{R}_{N,p}$.\\
The following lemmas give the rate of convergence of each of these terms.

\begin{lem} \label{l:varunblind}
The following upper bound holds:
$$\sup_{\onto{ Z_B   = a_B^TX_B}{ \left\|a_B\right\|_2 = 1}} \mathbb{E} \Bigg[      \Big( \bar{Z}_{B}  - \hat{Z}_{B} \Big)^2               \Bigg] ^{\frac{1}{2}}\leq  \frac{\sqrt{M}(m+M)}{m^2}  r_N. $$
\end{lem}

\begin{lem}\label{l:biais}
The following upper bound holds:
$$\sup_{ \onto{ Z_B   = a_B^TX_B}{ \left\|a_B\right\|_2 = 1}} \mathbb{E} \Bigg[      \Big( \bar{Z}_{B}  - \tilde{Z}_{B} \Big)^2               \Bigg] ^{\frac{1}{2}}\leq \frac{M}{m^4}  \left( \frac{M^{\frac{5}{2}}}{m} + M^{\frac{3}{2}} \right)   m_p. $$
\end{lem}

%
%

The proofs of these Lemmas are postponed to the Appendix

\end{proof}

\section{The blind case}

Assume now that there is only one sample $X$ available, observed on $O_N$.  In this case, the two previous estimation steps are linked and the asymptotics between the estimation of the spectral density and the projection step must be carefully chosen.\\ \indent In particular, we must pay a special attention to the behavior of the variance term to balance the two errors.  Indeed, to perform the prediction, we will use the whole available sample for the estimation, but then chose a window $O_{p(N)} \subset O_N$, for a suitable $p(N)$ to make the prediction. Doing this, we may go beyond the dependency problem induced by the fact that the same sample has to be used for both estimation and filtering.\vskip .1in

The following theorem provides the rate of convergence of the error term.

\begin{thm}
Assume that $\hat{f}_{p(N)}$ is built with the observation sample $X$ and that Assumptions [A1] and [A2] holds. The risk admits the following upper bound
$$\mathcal{R}_{N,p(N)} \leq   \frac{m+M}{m^2}  \mathbb{E} \left[ \left( \sum_{i \in O_{p(N)}} X_i^2 \right)^2        \right]^{\frac{1}{4}}   r_{N} +  \frac{M}{m^4}  \left( \frac{M^{\frac{5}{2}}}{m} + M^{\frac{3}{2}} \right)   m_{p(N)} .$$
\end{thm}
We point out that we obtain two terms alike those in Theorem~\ref{th:indep}. The main difference comes from the extra term $$\mathbb{E} \left[ \left( \sum_{i \in O_{p(N)}} X_i^2 \right)^2        \right]^{\frac{1}{4}}, $$
which corresponds to  the price to pay to use the same observation sample.  Using $\mathbb{E}[X_i^2X_j^2] \leq 3 \mathbb{E}[X_i^2]\mathbb{E}[X_j^2]$, leads to the following upper bound
\begin{align*} \mathbb{E} \left[ \left( \sum_{i \in O_{p(N)}} X_i^2 \right)^2        \right]^{\frac{1}{4}} & \leq 3 \sum_{i \in O_{p(N)}} \left(\mathbb{E} \left[   X_i^2 \right] \right)^2
\\ & \leq (\sharp O_{p(N)})^2 ., \end{align*} where $\sharp $ denotes the cardinal of a set. Finally, we see that the two error terms are linked in an opposite way such that the window parameter $p(N)$ must be chosen to balance the bias and the variance by minimizing with respect to $p(N)$ the quantity
\begin{equation} \label{eq:trade}
\frac{m+M}{m^2} \sqrt{\sharp O_{p(N)}}  r_{N} +   \frac{M}{m^4}  \left( \frac{M^{\frac{5}{2}}}{m} + M^{\frac{3}{2}} \right)   m_{p(N)} .\end{equation}

For instance, consider the special case of a process defined on $\mathbb{Z}^2$. Assume that $B$ is a finite subset of $\mathbb{Z}^2$, and that the we have at hand the following sequence of nested sub graphs $G_N = [-N,N]^2$. Then we get the following approximations for the quantities of Equation~\ref{eq:trade}, for some constants $C,k$, 
\begin{itemize}
 \item $\sharp O_N $ is of order $ 4N^2$
\item $r_N$ is of order $ \frac{C}{N}$ (see for instance \cite{dahlhaus})
\item $m_p$ is of order $ \frac{1}{p-k}$.
\end{itemize}

Then minimizing~\eqref{eq:trade} implies minimizing
$$C_1  \frac{p(N)}{N} + C_2\frac{1}{p(N)-k},  $$
which is achieved for $p(N) \approx \sqrt{N}$. Therefore the error is such that $$\mathcal{R}_{N,p(N)} = O\left( \frac{1}{\sqrt{N} }\right).$$ Here, the blind case leads with this method to an important loss since the error in the case of the independent sample would have been of order $\frac{1}{N}$.  A lower bound would be necessary to fully clarify this result, however obtaining such a bound seems a very difficult task which falls beyond the scope of this paper. Nevertheless, to our knowledge, the blind prediction of a graph indexed random process  has never been tackled before.

\begin{proof}

The only thing which remains to be calculated is given by the following lemma.

\begin{lem} \label{l:varblind}
The following upper bound holds:
$$\sup_{\onto{ Z_B   = a_B^TX_B}{ \left\|a_B\right\|_2 = 1}} \mathbb{E} \Bigg[      \Big( \bar{Z}_{B}  - \hat{Z}_{B} \Big)^2               \Bigg] ^{\frac{1}{2}}\leq   \frac{m+M}{m^2}  \mathbb{E} \left[ \left( \sum_{i \in O_{p(N)}} X_i^2 \right)^2        \right]^{\frac{1}{4}}   r_N. $$
\end{lem}
The proof is postponed in Appendix

\end{proof}

\section{Appendix} \label{s:append}

\begin{proof} of Lemma \ref{l:varunblind}

First, denote 
$$A_{N,p}\left(f,\hat{f}_N\right) =  \K_{BO_p}(\hat{f}_N) \left(\K_{O_p}(\hat{f}_N)\right)^{-1} - \K_{BO_p}(f) \left(\K_{O_p}(f)\right)^{-1} ,$$
and note that, $\mathbb{P}_X \otimes \mathbb{P}_Y$-a.s.,
$$\left(\hat{Z}_{B} - \bar{Z}_{B}\right)^2 = a_{B}^T\left(A_{N,p}\left(f,\hat{f}_N\right)\right)^TX_{O_p} X_{O_p}^TA_{N,p}\left(f,\hat{f}_N\right)a_{B}$$

Then, notice that 
\begin{align*}
 \sup_{\onto{ Z_B   = a_B^TX_B}{ \left\|a_B\right\|_2 = 1}} \mathbb{E} \Bigg[      \Big( \bar{Z}_{B} & - \hat{Z}_{B} \Big)^2               \Bigg] ^{\frac{1}{2}}
 \\ &  \leq   \left( \sup_{ \onto{ \operatorname{supp}(a_{B}) \subset{B}}  { \left\| a_{B} \right\|_2 = 1}} \mathbb{E} _{\mathbb{P}_Y}\left[  a_{B}^T \left(A_{N,p}\left(f,\hat{f}_p\right)\right)^T\K_{O_p}(f)A_{N,p}\left(f,\hat{f}_N\right) a_{B}  \right]           \right)^{\frac{1}{2}} 
\\ & \leq  \mathbb{E} _{\mathbb{P}_Y}\left[ \sup_{ \onto{ \operatorname{supp}(a_{B}) \subset{B}}  { \left\| a_{B} \right\|_2 = 1}}  a_{B}^T    \left(A_{N,p}\left(f,\hat{f}_N\right)\right)^T\K_{O_p}(f)A_{N,p}\left(f,\hat{f}_N\right) a_{B}  \right]  ^{\frac{1}{2}} .
 \end{align*}

Thus,

$$  \sup_{ \onto{ Z_B   = a_B^TX_B}{ \left\|a_B\right\|_2 = 1}} \mathbb{E} \left[      \left( \bar{Z}_{B} - \hat{Z}_{B} \right)^2               \right] ^{\frac{1}{2}} \leq 
\mathbb{E} _{\mathbb{P}_Y}\left[
  \left\| \left( A_{N,p}\left(f,\hat{f}_N\right)   \right)^T\K_{O_p}(f)A_{N,p}\left(f,\hat{f}_N\right)            \right\|_{2,op} \right]^{\frac{1}{2}}.$$

Furthermore, it holds, $\mathbb{P}_Y$-a.s., that
\begin{align*}
  \left\| \left( A_{N,p}\left(f,\hat{f}_N\right)   \right)^T\K_{O_p}(f)A_{N,p}\left(f,\hat{f}_N\right)            \right\|_{2,op}  & \leq 
    \left\|  A_{N,p}\left(f,\hat{f}_N\right)   \right\|_{2,op}^2
      \left\| \K_{O_p}(f)            \right\|_{2,op} \end{align*}

      But,
      \begin{align*}
    \left\|  A_{N,p}\left(f,\hat{f}_N\right)   \right\|_{2,op}     
        & \leq          \left\|         \K_{BO_N}(\hat{f}_N) \left(\K_{O_p}(\hat{f}_N)\right)^{-1}  -  \K_{BO_p}(f) \left(\K_{O_p}(\hat{f}_N)\right)^{-1} \right\|_{2,op}\\ & +    \left\| \K_{BO_p}(f) \left(\K_{O_p}(\hat{f}_N)\right)^{-1} - \K_{BO_p}(f) \left(\K_{O_p}(f)\right)^{-1} \right\|_{2,op}
\\ & \leq   \left\|         \left(\K_{O_p}(\hat{f}_N)\right)^{-1}  \right\|_{2,op}    \left\|         \K_{BO_p}(\hat{f}_N)  -  \K_{BO_p}(f) \right\|_{2,op} + 
  \left\| \K_{BO_p}(f)   \right\|_{2,op}    \\ & \times       \left\|        \left(\K_{O_p}(f)\right)^{-1} \right\|_{2,op}   \left\|        \left(\K_{O_p}(\hat{f}_N)\right)^{-1} \right\|_{2,op} 
 \left\|  \K_{O_p}(\hat{f}_N)- \K_{O_p}(f) \right\|_{2,op}
      \\ & \leq \frac{1}{m} \left\| f - \hat{f}_N \right\|_\infty + \frac{M}{m^2}  \left\| f - \hat{f}_N \right\|_\infty 
        \end{align*}
Here we used the inequality
$$\left\| \mathcal{K}(f) \right\|_{2,op} \leq \left\| f \right\|_\infty $$

We get 
\begin{align*}
 \sup_{ \onto{ Z_B   = a_B^TX_B}{ \left\|a_B\right\|_2 = 1}} \mathbb{E} \left[      \left( \bar{Z}_{B} - \hat{Z}_{B} \right)^2               \right] ^{\frac{1}{2}} & \leq \frac{\sqrt{M}(m+M)}{m^2}  
\mathbb{E} _{\mathbb{P}_Y}\left[
                                                                     \left\| f - \hat{f}_N \right\|_\infty ^2      \right]^{\frac{1}{2}}.
                                                                     \\ & \leq \frac{\sqrt{M}(m+M)}{m^2}  r_N
\end{align*}
\end{proof}

\begin{proof} of Lemma \ref{l:biais}

First, define for all $A \subset G$, the operator $p_{A}$ by 
$$ \forall i,j \in G, \left(p_{A}\right)_{ij} = \ind_{i \in A}\ind_{i = j}.$$
To compute the rate of convergence of the bias term $\bar{Z}_B-\tilde{Z}_B$, we can compute directly, 


\begin{align*}
 \sup_{ \onto{ Z_B   = a_B^TX_B}{ \left\|a_B\right\|_2 = 1}} \mathbb{E} \Bigg[      \big( \bar{Z}_{B}& - \tilde{Z}_{B} \Big)^2               \Bigg] ^{\frac{1}{2}} 
 \\ & \leq \left\|    \K_{BO_p}(f) \left(\K_{O_p}(f)\right)^{-1} p_{O_p}  -   K_{BO}(f) \left(\K_{O}(f)\right)^{-1}             \right\|_{2,op} \left\| \K_{O}(f)\right\|_{2,op}^{\frac{1}{2}}
\end{align*}

Note that, since $B \cup O = G$, we have immediately that $\left( \K_{B \cup O}(f) \right)^{-1} = \K(\frac{1}{f})$.

Using a Schur decomposition, we get
\begin{align*}
 \sup_{\onto{ Z_B   = a_B^TX_B}{ \left\|a_B\right\|_2 = 1}} \mathbb{E} \Bigg[      \big( \bar{Z}_{B}& - \tilde{Z}_{B} \Big)^2               \Bigg] ^{\frac{1}{2}} 
 \\ & \leq \sqrt{M} \left\|     \K_{BO_p}(f) \left(\K_{O_p}(f)\right)^{-1}   +   \left( \K_{B}(\frac{1}{f}) \right)^{-1} \K_{BO}(\frac{1}{f})                     \right\|_{2,op} 
 \\
  \\ & \leq \sqrt{M} \left\|     \K_{BO_p}(f) \left(\K_{O_p}(f)\right)^{-1}   +  \left(  \K_{B}(\frac{1}{f})  \right)^{-1}\K_{BO_p}(\frac{1}{f})                \right\|_{2,op} 
 \\ & +  \sqrt{M}    \left\| \left(  \K_{B}(\frac{1}{f})  \right)^{-1}  \right\|_{2,op}  \left\|     \K_{B(O\setminus O_p)}(\frac{1}{f})                 \right\|_{2,op}
 \\
 \\ & \leq \sqrt{M} \left\|     \K_{B}(\frac{1}{f}) \K_{BO_p}(f)   +  \K_{BO_N}(\frac{1}{f})   \K_{O_p}(f)          \right\|_{2,op}
 \\ & \times   \left\|     \left(\K_{O_p}(f)\right)^{-1}                 \right\|_{2,op}   \left\|    \left( \K_{B}(\frac{1}{f}) \right)^{-1}              \right\|_{2,op}  + 
 M^{\frac{3}{2}} \left\|   \K_{B(O\setminus O_p)}(\frac{1}{f})                    \right\|_{2,op}
 \\
   \\ & \leq \frac{M^{\frac{3}{2}}}{m} \left\|   - \K_{B(O\setminus O_p)}(\frac{1}{f})   \K_{(O\setminus O_p)O_p}(f)      \right\|_{2,op}
 \\ &   + 
 M^{\frac{3}{2}} \left\|    \K_{B(O\setminus O_p)}(\frac{1}{f})                  \right\|_{2,op}
 \end{align*}

 This leads to 
 
 \begin{align}
  \sup_{ \onto{ Z_B   = a_B^TX_B}{ \left\|a_B\right\|_2 = 1}} \mathbb{E} \Bigg[      \big( \bar{Z}_{B}& - \tilde{Z}_{B} \Big)^2               \Bigg] ^{\frac{1}{2}}
   \leq \left( \frac{M^{\frac{5}{2}}}{m} + M^{\frac{3}{2}} \right) \left\|     \K_{B(O \setminus O_p)}(\frac{1}{f})                     \right\|_{2,op}
%
%
 \end{align}

Moreover, we can write 

\begin{align*}
\left\|     \K_{B(O \setminus O_p)}(\frac{1}{f})         \right\|_{2,op} &\leq \sum_{k \geq 0} \left| (\frac{1}{f})_k \right| \left\| \left(A^k\right)_{B(O \setminus O_p)}   \right\|_{2,op}
\\ & \leq  \sum_{k \geq d_G\left(B, (O\setminus O_p)\right)} \left|  (\frac{1}{f})_k  \right| 
\\ & \leq  \sum_{k \geq d_G\left(B, (G\setminus G_p)\right)} \left|  (\frac{1}{f})_k  \right| 
\\ & \leq  \frac{1}{d_G\left(B, (G\setminus G_p)\right)} \sum_{k \geq d_G\left(B, (G\setminus G_p)\right)} k \left|  (\frac{1}{f})_k  \right| 
\\ & \leq \frac{M}{m^2}m_p
\end{align*}

\end{proof}

\begin{proof} of Lemma \ref{l:varblind}

For sake of simplicity, let us denote $O_p$ instead of $O_{p(N)}$ in the whole proof.

We still denote 
$$A_{N,p}\left(f,\hat{f}_N\right) =  \K_{BO_p}(\hat{f}_N) \left(\K_{O_p}(\hat{f}_N)\right)^{-1} - \K_{BO_p}(f) \left(\K_{O_p}(f)\right)^{-1} ,$$

Then, notice that 
\begin{align*}
 \sup_{\onto{ Z_B   = a_B^TX_B}{ \left\|a_B\right\|_2 = 1}} &\mathbb{E} \Bigg[      \Big( \bar{Z}_{B}  - \hat{Z}_{B} \Big)^2               \Bigg] ^{\frac{1}{2}}
 \\ &  \leq   \left( \sup_{ \onto{ \operatorname{supp}(a_{B}) \subset{B}}  { \left\| a_{B} \right\|_2 = 1}} \mathbb{E} \left[\operatorname{Tr} \left(  a_{B}^T \left(A_{N,p}\left(f,\hat{f}_N\right)\right)^TX_{O_p}X_{O_p}^TA_{N}\left(f,\hat{f}_N\right) a_{B} \right)  \right]           \right)^{\frac{1}{2}} 
 \\ &  \leq   \left( \sup_{ \onto{ \operatorname{supp}(a_{B}) \subset{B}}  { \left\| a_{B} \right\|_2 = 1}} \mathbb{E} \left[\operatorname{Tr} \left(  X_{O_p}X_{O_p}^TA_{N,p}\left(f,\hat{f}_N\right) a_{B} a_{B}^T \left(A_{N,p}\left(f,\hat{f}_N\right)\right)^T\right)  \right]           \right)^{\frac{1}{2}}
 \end{align*}
 
 Applying Cauchy Schwartz inequality, we get
 \begin{align*}
  \sup_{\onto{ Z_B   = a_B^TX_B}{ \left\|a_B\right\|_2 = 1}} \mathbb{E} \Bigg[      \Big( \bar{Z}_{B} & - \hat{Z}_{B} \Big)^2               \Bigg] ^{\frac{1}{2}}
\\  & \leq    \Bigg(        \mathbb{E} \left[\operatorname{Tr} \left( \left( X_{O_p}X_{O_p}^T \right)^2      \right)  \right]  
\\ & \times             \sup_{ \onto{ \operatorname{supp}(a_{B}) \subset{B}}  { \left\| a_{B} \right\|_2 = 1}} \mathbb{E} \left[\operatorname{Tr} \left( \left(A_{N,p}\left(f,\hat{f}_N\right) a_{B} a_{B}^T \left(A_{N,p}\left(f,\hat{f}_N\right)\right)^T\right)^2 \right)  \right]           \Bigg)^{\frac{1}{4}}
 \end{align*}

But, on the one hand, we have
\begin{align*}
\mathbb{E} \left[\operatorname{Tr} \left( \left( X_{O_p}X_{O_p}^T \right)^2      \right)  \right]  & = \mathbb{E} \left[\operatorname{Tr} \left( \left( X_{O_p}^TX_{O_p} \right)^2      \right)  \right]  
\\ & =  \mathbb{E} \left[ \left( \sum_{i \in O_p} X_i^2 \right)^2        \right]  
\end{align*}
 
 And on the other hand, we can write
 
\begin{align*}
  \sup_{ \onto{ \operatorname{supp}(a_{B}) \subset{B}}  { \left\| a_{B} \right\|_2 = 1}} & \mathbb{E} \left[\operatorname{Tr} \left( \left(A_{N}\left(f,\hat{f}_N\right) a_{B} a_{B}^T \left(A_{N,p}\left(f,\hat{f}_N\right)\right)^T\right)^2 \right)  \right]         
\\ & \leq   \sup_{ \onto{ \operatorname{supp}(a_{B}, b_{B}) \subset{B}}  { \left\| a_{B} \right\|_2 = \left\| b_{B} \right\|_2 = 1}} \mathbb{E} \Bigg[\operatorname{Tr} \bigg(  a_{B}^T \left(A_{N,p}\left(f,\hat{f}_N\right)\right)^TA_{N,p}\left(f,\hat{f}_N\right)
\\ &\times b_{B} b_{B}^T  \left(A_{N,p}\left(f,\hat{f}_N\right)\right)^TA_{N,p}\left(f,\hat{f}_N\right)a_{B}   \bigg)  \Bigg]         
\\ & \leq   \sup_{ \onto{ \operatorname{supp}(b_{B}) \subset{B}}  {  \left\| b_{B} \right\|_2 = 1}}  \mathbb{E} \Bigg[       \left\| A_{N,p}\left(f,\hat{f}_N\right) \right\|_{2,op}^4  \left\|b_{B} b_{B}^T  \right\|_{2,op}                                                \Bigg]       
\\ & \leq   \mathbb{E} \Bigg[       \left\| A_{N,p}\left(f,\hat{f}_N\right) \right\|_{2,op}^4                                       \Bigg]  
\\ & \leq   \left( \frac{m+M}{m^2}\right)^4  \mathbb{E} \Bigg[       \left\| f -\hat{f}_N  \right\|_{\infty}^4                                       \Bigg]  
\end{align*}

Thus,

\begin{align*}
  \sup_{ \onto{ Z_B   = a_B^TX_B}{ \left\|a_B\right\|_2 = 1}} \mathbb{E} \Bigg[      \Big( \bar{Z}_{B} & - \hat{Z}_{B} \Big)^2               \Bigg] ^{\frac{1}{2}}
 \leq \frac{m+M}{m^2}  \mathbb{E} \left[ \left( \sum_{i \in O_p} X_i^2 \right)^2        \right]^{\frac{1}{4}}   r_N
 \end{align*}

 \end{proof}

\bibliographystyle{abbrv}
\bibliography{biblioGraphKrig}

\end{document}